\magnification=\magstep1
\advance\voffset by 1 true cm
\baselineskip=15pt
\overfullrule =0pt 
\font\bbb=msbm10
\def\R{\hbox{\bbb R}}   
\def\N{\hbox{\bbb N}}    
\def\C{\hbox{\bbb C}}
\def\P{\hbox{\bbb P}}

\centerline {\bf ARRANGEMENTS, MILNOR FIBERS and POLAR CURVES}

\vskip1.5truecm
\centerline {\bf by Alexandru Dimca }

\vskip1.5truecm

{\bf 1. The main results}

\bigskip

Let ${\cal A}$ be a hyperplane arrangement in the complex projective space $\P^n$, with $n>0$. Let $d>0$ be the number of hyperplanes in this arrangement and choose a linear equation $H_i: \ell_i(x)=0$ for each hyperplane $H_i$ in ${\cal A}$,
for $i=1,...,d.$

Consider the homogeneous polynomial $Q(x)= \prod _{i=1,d}\ell_i(x) \in \C[x_0,...,x_n]$ and the corresponding principal open set $D(Q)=\P^n \setminus \cup _{i=1,d}H_i$. The topology of the hyperplane arrangement complement $D(Q)$ is a central object of study in the theory of hyperplane arrangements, see Orlik-Terao [OT1].

There is a gradient map associated to any homogeneous polynomial $h \in \C[x_0,...,x_n]$, namely

$$  grad(h): D(h) \to \P^n, ~~~~(x_0:...:x_n) \mapsto (h_0(x):...:h_n(x))$$
where $D(h)= \{x \in \P^n ; h(x) \not= 0\}$ and $h_i={\partial h \over \partial x_i}$. A nice consequence of our main result is the following.

\bigskip

{\bf Theorem 1.}  {\it For any projective arrangement ${\cal A}$  as above one has
$$b_n(D(Q))=deg(grad(Q)).$$
In particular, the following are equivalent.

(i) the morphism $grad(Q)$ is dominant;

(ii)  $b_n(D(Q)) >0$ and

(iii) the projective arrangement ${\cal A}$ is essential, i.e. the intersection $\cap _{i=1,d}H_i$ is empty.}

\bigskip
Our main result stated as Theorem 2. below is a substantial improvement of some of  the main results by Orlik and Terao [OT2]. Let ${\cal A}'$ be the affine hyperplane arrangement in $\C^{n+1}$ associated to the projective arrangement ${\cal A}$. Note that $Q(x)=0$ is a reduced equation for the union $N$ of all the hyperplanes in ${\cal A}'$.

Let $f \in \C[x_0,...,x_n]$ be a homogeneous polynomial of degree $e>0$ with
global Milnor fiber $F=\{ x \in \C^{n+1}| f(x))=1\}$, see for instance [D] for more on such varieties. Let $g: F \setminus N \to \R$
be the function $g(x)=Q(x) {\overline Q}(x)$. The polynomial $f$ is called 
${\cal A}'$-generic if

(GEN1) the restriction of $f$ to any intersection $L$ of hyperplanes in ${\cal A}'$ is non-degenerate, in the sense that the associated projective hypersurface in $\P(L)$ is smooth, and

(GEN2) the function $g$ is a Morse function.

Orlik and Terao have shown in [OT2] that for an {\it essential} arrangement ${\cal A}'$, the set of ${\cal A}'$-generic functions $f$ is dense in the set of homogeneous polynomials of degree $e$, and, as soon as we have an ${\cal A}'$-generic function $f$, the following basic properties hold for any arrangement.

(P1) $~~b_q(F,F \cap N)=0$ for $q \not= n$ and

(P2) $~~b_{n}(F,F \cap N) \leq |C(g)|$, where $C(g)$ is the critical set of the Morse function $g$.

Moreover, for a special class of arrangements called {\it pure} arrangements it is shown in [OT2] that (P2) is actually an equality.
Note that (P1) and basic additivity properties of Euler characteristics, see for instance [DL], imply that (P2) is equivalent to

(P2') $~~~~~~~(-1)^n \chi(F\setminus N) \leq |C(g)|$.

With this notation our main result is the following.

\bigskip

{\bf Theorem 2.}

{\it For  any arrangement ${\cal A}'$ the following hold

(i) The set of ${\cal A}'$-generic functions $f$ is dense in the set of homogeneous polynomials of degree $e>1$;

(ii)  $~~~~\chi(F \setminus N) = (-1)^n|C(g)|$.}

\bigskip

In fact, Theorem 1. follows from Theorem 2. by taking $f$ a general linear form.
Our method of proof for Theorem 2. is completely different from the methods used by Orlik-Terao (though in both cases  there is  some Morse theory in the background) and uses
 the affine Lefschetz theory developped by N\'emethi in [N1-2].

In section 2. we recall the main results of [N1-2], emphasizing  their reformulations in   [CD] in terms of polar curves and, in order to prepare the reader for the more complicated proof of Theorem 2. given in section 3., we prove the following
  topological description for the degree of the gradient $grad(h)$ of any homogeneous polynomial $h$.

\bigskip

{\bf Proposition 3.} {\it For any homogeneous polynomial   $h \in \C[x_0,...,x_n]$, one has

$deg(grad(h))= (-1)^n \chi (D(h) \setminus H)$, with $H$ a general hyperplane
in $\P^n$. }

\bigskip
Via a simple  (surely known) fact  concerning the topology of hyperplane arrangements, see Lemma 7., Proposition 3. gives another, more direct proof of Theorem 1.
\bigskip
In the last section we clarify the proof of a recent result by R. Randell concerning the minimality of the complement $D(Q)$, a basic result in view of the applications, see [PS], [R].

\bigskip

The author thanks Stefan Papadima for raising the question answered by Theorem 1 above
and for lots of helpful comments.  In particular he informed me that Theorem 1 was proved by Paltin Ionescu in the case $n=2$ by completely different methods.

He also thanks Pierrette Cassou-Nogu\`es for drawing his attention on R. Randell's preprint.

\bigskip
{\bf 2. Polar curves and degree of gradient maps}

\bigskip

The use of the local polar varieties in the study of singular spaces
is already a classical subject, see L\^e-Teissier [LT] and the
references therein. If fact, all the results in this section can be obtained from the local results of L\^e [L\^e], but we prefer to devellop the general theory which is needed in section 3.

Global polar curves in the study of the topology of polynomials is a topic under intense investigations, see for instance Cassou-Nogu\`es and Dimca [CD], Siersma and Tib\u ar [ST], [T].

We recall briefly the notation and the results from [CD]. Let $h \in \C[x_0, ... ,x_n]$ be a polynomial (even non-homogeneous) and assume that the
fiber $F_t=h^{-1}(t)$ is smooth and connected, for some fixed $t \in \C$.

 For any hyperplane in $\P^n$, $ H:\ell =0$ where $  \ell(x)=h_0x_0+h_1x_1+...+h_nx_n $
we define the corresponding polar variety $\Gamma _H $ to be the union
of the irreducible components of the variety
$$ \{x \in \C^{n+1} \ |  \  \  rank(dh(x),d\ell (x))=1 \} $$
which are not contained in the critical set $S(h) = \{x \in \C^{n+1} \ | \
dh(x)=0 \} $ of $h$.

\bigskip

\noindent {\bf Lemma 4.} (see [CD], [ST])

{\it For a generic hyperplane $H$ we have the following properties.

(i) The polar variety $\Gamma _H$ is either empty or a curve,
i.e. each
irreducible component of $\Gamma _H$ has dimension 1.

(ii) dim$(F_t \cap \Gamma _H) \leq 0$ and the intersection
multiplicity
$(F_t,\Gamma_H)$ is independent of $H$.

(iii) The multiplicity $(F_t, \Gamma _H)$ is equal to the number of
tangent hyperplanes to $F_t$ parallel to the hyperplane $H$. For each
such tangent hyperplane $H_a$, the intersection $ F_t \cap H_a$ has
precisely one singularity, which is an ordinary double point.}

\bigskip

\noindent {\bf Definition 5.}

{\it The non-negative integer $(F_t, \Gamma _H)$ is called the polar
invariant
of the hypersurface $F_t$ 
and is denoted by $P(F_t)$.}
\bigskip
Note that $P(F_t)$ corresponds exactly to the classical notion of
class of a projective hypersurface, see [L].

We think of a projective hyperplane $H$ as above as the direction of an affine hyperplane $H' =\{ x \in \C^{n+1}|  \ell (x)=s \}$ for $s \in \C$. All the hyperplanes with the same direction form a pencil, and it is precisely the pencils of this type that are used in the affine Lefschetz theory, see [N1-2].
One of the  main results  in [CD] is the following, see also [ST] or [T] for similar results.
\bigskip

\noindent{\bf Proposition 6.}

\bigskip

{\it For a generic hyperplane $H'$ in the pencil of all hyperplanes in $\C^{n+1}$ with a fixed generic direction $H$, the homotopy type of the fiber $F_t$ is
obtained from the homotopy type of the section $F_t \cap H'$ by
attaching $P(F_t)$  cells of dimension $n$.}

In particular
$$P(F_t)=(-1)^{n}(\chi(F_t)- \chi(F_t \cap H'))=(-1)^n \chi(F_t \setminus H')$$
Moreover in this statement 'generic' means that the hyperplane $H'$ has to verify the following two conditions.

(g1) its direction, which is the hyperplane in $\P^n$ given by the homogeneous part of degree one in an equation for $H'$ has to be generic, and

(g2) the intersection $ F_t \cap H'$ has to be smooth.

These two conditions are not stated in [CD], but the reader should have no problem in checking them by using Theorem 3' in [CD] and the fact proved by N\'emethi in [N1-2] that the only bad sections in a good (i.e. the analog of a Lefchetz pencil in the projective Lefschetz theory, see [L]) pencil are the singular sections. Completely similar results hold for generic pencils with respect to a closed smooth subvariety $Y$ in some affine space $\C^N$, see [N1-2], but note that the polar curves are not mentionned there. In the next section we will need this more general setting.

Assume from now on that the polynomial $h$ is homogeneous of degree $d$ and that $t \not=0$. It follows from (g1) and (g2) above that we may choose the generic hyperplane $H'$ passing through the origin.

Moreover, in this case, the polar curve $\Gamma _H$, being defined by homogeneous equations, is a union of lines $L_j$ passing through the origin. For each such line we choose a parametrization $t \mapsto a_j t$ for some $a_j \in \C^{n+1}, a_j \not=0$. It is easy to see that the intersection $F_t \cap L_j$ is either empty (if $h(a_j)=0$) or consists of exactly $d$ distinct points with multiplicity one (if $h(a_j)\not=0$). The lines of the second type are in bijection with the points in $grad(h)^{-1}(D_{H'})$,
where $D_{H'}\in \P^n$ is the point corresponding to the direction of the hyperplane $H'$. It follows that
$$ d \cdot deg(grad(h))= P(F_t).$$

The $d$-sheeted unramified coverings $ F_t \to D(h)$ and $F_t \cap H' \to D(h) \cap H$ give the result, where $H$ is the projective hyperplane corresponding to the affine hyperplane (passing through the origin) $H'$. Indeed, they imply the equalities: $\chi(F_t)=d \cdot \chi (D(h))$ and $\chi(F_t \cap H')=d \cdot \chi (D(h) \cap H)$.

\bigskip

{\bf Lemma 7.} {\it  For any arrangement ${\cal A}$  as above one has $(-1)^n
\chi (D(f) \setminus H)=b_n(D(f))$.}
\bigskip

Proof.

Here we just give the main idea, since the details are standard. One has to use the method of deletion and restriction, see [OT1], p. 17, the obvious additivity of the Euler characteristics  and, more subtly, the additivity of the top Betti numbers coming from the exact sequence (8) in [OT1], p. 20 or (3.8) in [DL].

\bigskip

To complete the proof of Theorem 1 we still have to explain why the claims $(ii)$ and $(iii)$ are equivalent. If the projective arrangement is not essential, then using a projection onto $\P^{n-1}$ with center a point in all the hyperplanes $H_i$ we get a fiber bundle $D(Q) \to U$ with fiber $\C$ and base $U$, an affine variety of dimension $n-1$. This implies $b_n(D(Q))=0$.

If the arrangement is essential, then $d \geq n+1$ and we may assume that
$\ell _i(x)=x_{i-1}$ for $i=1,...,n+1.$ In the case $d = n+1$, we are done, since in this case $D(Q)= (\C^*)^n$ and hence $b_n(D(Q))=1$. In the remaining case 
$d>n+1$, one should use the additivity of the top Betti numbers alluded above.
\bigskip

{\bf 3. Proof of Theorem 2.}

\bigskip

Let $f \in \C[x_0,...,x_n]$ be for the moment any non-zero homogeneous polynomial of degree $e$. Note that the $e$-sheeted covering $F \to D(f)$ induces an $e$-sheeted covering $F \setminus N \to D(fQ)$. In particular
$$\leqno (8)~~~~~~~~~~ \chi (F \setminus N)=e \cdot \chi (D(fQ)).$$
In a similar way, if we set $X=Q^{-1}(1)$ and $X_0=\{x \in X|f(x)=0\}$ then there is a $d$-sheeted covering $X \setminus X_0 \to D(fQ)$ giving
$$\leqno (9)~~~~~~~~~~ \chi (X \setminus X_0)=d \cdot \chi (D(fQ)).$$
Let  $v: \C^{n+1} \to \C^N$ be the Veronese mapping of degree $e$ sending $x$ to all the monomials of degree $e$ in $x$ and set $Y=v(X)$.
Then $Y$ is a smooth closed subvariety in $\C^N$ and $v:X \to Y$ is an unramified (even Galois) covering of degree $c$, where $c=g.c.d.(d,e)$. To see this, use the fact that $v$ is a closed immersion on $\C^N \setminus \{0\}$ and $v(x)=v(x')$ iff $x'=u \cdot x$ with $u^c=1$.

Let $H$ be a generic hyperplane direction in $\C^N$ with respect to the subvariety $Y$ and let $C(H)$ be  the finite set of all the points $p \in Y$ such that there is an affine hyperplane $H'_p$ in the pencil determined by $H$ that is tangent to $Y$ at the point $p$ and the intersection $Y \cap H'_p$ has a complex Morse (alias non-degenerated, alias $A_1$) singularity. We can also assume that the affine hyperplane $H_0$ in this pencil is generic (i.e. $Y_0=Y \cap H_0$ is smooth) and then by [N1-2] we get
$$ \leqno (10) ~~~~~ \chi (Y,Y_0)=(-1)^n|C(H)|. $$
Under the Veronese mapping $v$, the generic hyperplane direction $H$ corresponds to a homogeneous polynomial of degree $e$ which we call from now on $f$.

To prove the first claim (i) note that the first condition (GEN1) is clearly generic, while the second condition (GEN2) is fulfilled by our polynomial $f$ above. Indeed, in view of the last statement at the end of the proof of Lemma (2.5) in [OT2] $g$ is a Morse function iff each critical point of $Q:F \setminus N \to \C$ is an $A_1$-singularity. Using the homogeneity of both $f$ and $Q$, this last condition on $Q$ is equivalent to the fact that each critical point of the function $f:X \to \C$ is an $A_1$ singularity, condition fulfilled in view of the choice of $H$ and since $v:X \to Y$ is a local isomorphism.

Now we pass on to the proof of the claim (ii) in Theorem 2. Under the Veronese mapping $v$, the pencil of hyperplanes in $\C^N$ gives rise to a pencil of hypersurfaces $V_s: f(x)=s$ in $\C^{n+1}$ such that $X_0=X \cap V_0$.
Let $C(V)$ be the finite set of all the points $q \in X$ such that there is an affine hypersurface $V_s$  that is tangent to $X$ at the point $q$. The $c$-sheeted covering $X \to Y$ induces a covering $C(V) \to C(H)$, and hence multiplying (10) by $c$ we get
$$ \leqno (11) ~~~~~ \chi (X,X_0)=(-1)^n|C(V)|. $$
Lemma (2.3) in [OT2] gives the following description of the critical set $C(g)$.
$$ \leqno (12) ~~~~~ C(g)=\{x \in F \setminus N| \ \  rank(dQ(x),df(x)) =1 \}. $$
To compare $|C(g)|$ and $|C(V)$ we proceed as follows. Let 
 the  polar variety $\Gamma _{Q,f} $ of the pair $(Q,f)$  be the union
of the irreducible components of the variety
$$ \{x \in \C^{n+1} \ |  \  \  rank(dQ(x),df(x))=1 \} $$
which are not contained in the critical set $S(Q) = \{x \in \C^{n+1} \ | \
dQ(x)=0 \} $ of $Q$. Then, exactly as in the simpler situation described in section 2., $\Gamma _{Q,f} $ is a union of lines  $L_j$ passing through the origin. For each such line we choose a parametrization $t \mapsto a_j t$ for some $a_j \in \C^{n+1}, a_j \not=0$. It is easy to see that the intersection $X \cap L_j$ is non-empty  iff the intersection $(F\setminus N) \cap L_j$ is non- empty  and, if this is the case the first (resp. the second) intersection consists of $d$ (resp. $e$) points. Therefore

$$|C(g)|={e \over d}\cdot |C(V)|={e \over d}(-1)^n \chi (X \setminus X_0)=(-1)^ne\cdot  \chi(D(fQ))=(-1)^n\chi (F \setminus N).$$
\bigskip
{\bf Remark 8.}

Both Theorem 1 and Theorem 2 above remain true (with the same proof) when we replace the polynomial $Q$ by a more general polynomial $Q_{\bf m} (x)=\ell_1(x)^{m_1} \cdot ... \cdot \ell_d(x)^{m_d}$ for any ${\bf m}=(m_1,...,m_d) \in (\N^*)^d.$ Note that $D(Q)=D(Q_{\bf m})$, hence the integers $deg(grad(Q_{\bf m}))$ and $|C(g_{\bf m})|$, with $g_m(x)=Q_{\bf m}(x) {\overline Q}_{\bf m}(x)$  are independent of $\bf m$, a result similar to the results in [OT3].

\bigskip

{\bf 4. Minimality of hyperplane arrangements}

\bigskip

In this section we discuss a recent result by Randell [R], closely related to the topic considered so far. Let $M^*=D(Q)$ be the complement of the projective arrangement in $\P^n$ and let $M=M({\cal A'})$ be the complement of the corresponding affine central arrangmemt in $\C^{n+1}$.

We say that a topological space $Z$ is minimal if $Z$ has the homotopy type of a CW-complex $K$ whose number of $k$-cells equals $b_k(K)$ for all $k \in \N$.

The importance of this notion for the topology of hyperplane arrangements was recently discovered by S. Papadima and A. Suciu, see [PS] but also [R] for various applications.
In view of this, the following result is crucial.

\bigskip

{\bf Theorem 9.} {\it The complements $M$ and $M^*$ are minimal spaces.}

\bigskip

Proof. The proof of this result given by Randell  in [R] is as follows.

It is clearly enough to treat the projective complement $M^*$.
Then one notices that the Milnor fiber $X$  of the arrangement (denoted by $F$ in [R]) has the homotopy type of a space obtained from $X \cap H'$ by attaching $n$-cells, where $H':\ell (x)=0$ is a generic hyperplane in $\C^{n+1}$ passing through the origin. Randell obtains this result by using the Morse function $m(x)=|\ell (x) |$, following the local results by L\^e [L\^e]. This is indeed possible, since in a homogeneous situation there is no difference between the local and the global case.

Then Randell claims that  the corresponding CW-structure constructed inductively on the Milnor fiber $X$ is invariant with respect to the covering transformations of the projection $p:X \to M^*$ and hence it gives rise by taking quotients to a CW-structure on $M^*$. This is claimed in the Introduction, then in Example 3 and it is used in the proof of the main result (Theorem 4), but in my opinion no clear argument is given to support this claim.

In the following we propose two ways to avoid this claim (apparently difficult to prove): the new idea is to  construct directly the CW-structure on the base space $M^*$ without looking first for the CW-structures on the Milnor fiber $X$.

The first method is the simplest: using the Affine Lefschetz Theorem of Hamm, see Theorem 5 in [H], we know that for a generic projective hyperplane $H$, the space $M^*$ has the homotopy type of a space obtained from $M^* \cap H$ by attaching $n$-cells. The number of these cells is given by
$$(-1)^n \chi (M^*, M^* \cap H)=(-1)^n \chi (M^* \setminus H)=b_n(M^*)$$
see Lemma 7 above.

To finish the proof of the minimality of $M^*$ we proceed by induction using the equalities
$$b_k(M^*)=b_k(M^* \cap H)$$
for $0 \leq k <n$ which are easily proved as in Randell [R].

The second method to describe the CW-structure on the base space $M^*$ is more complicated but more precise (it is in fact the proof of the result by Hamm used above). Note first that the Morse function $m=|\ell |$ of $X$ relative to $X \cap H'$ is not proper, hence to get the results we have in fact to cut everything with a big closed ball $B_R=\{x \in \C^{n+1}; |x| \leq R\}$ and use Morse theory as in Hamm [H].

The group $G$ of covering transformation of the projection $p:X \to M^*$ is spanned by $T(x)=(ux_0,...,ux_n)$ with $u=exp(2 \pi i/d)$. It follows that the function $m$, being invariant under the group $G$, gives rise to a Morse function $m^*:M^* \to \R$ of $M^*$ relative to $M^* \cap H$, with $H$ the projective hyperplane associated to $H'$. Note that we are again in a non-compact situation, so to make the proof complete we have to replace $M^*$ by the manifold with boundary $p(X \cap B_R)$ for $R>>0$ and check that there are no singularities on the boundary.

This Morse function has  $b_n(M^*)$ critical points by Lemma 7 and all of them have index $n$ since $p$ is a local diffeomorphism. The proof is finished in the same way as in the first method.

\bigskip

{\bf Remark 10.}

We can use the second method to obtained an explicit Morse function $g^*$ on the open set $D(f)$ relative to the intersection
$D(f) \cap \cup_iH_i$=$D(f) \cap p(N)$, or an explicit Morse function  on $M^*$ relative to the intersection with a generic hypersurface
$M^* \cap \{x \in \P^n; f(x)=0\}$. However in the case $d=deg(f)>1$ we have $b_n(M^*) \not= (-1)^n \chi (M^* \setminus \{x \in \P^n; f(x)=0\})$ and hence this Morse function cannot be used to prove the minimality of $M^*$. It is also known that  the open set $D(f)$ is not minimal for $f$ generic.

\bigskip

\noindent {\bf REFERENCES}
\bigskip

\item{[CD]} Pi. Cassou-Nogu\`es, A. Dimca: Topology of complex polynomials via polar curves, Kodai Math. J. 22(1999), 131-139.

\item{[D]} A. Dimca: Singularities and Topology of Hypersurfaces,
Universitext,Springer, 1992.

\item {[DL]} A. Dimca, G.I. Lehrer: Purity and equivariant weight
polynomials, dans le volume: Algebraic Groups and Lie Groups, editor
G.I. Lehrer, Cambridge University Press, 1997.

\item{[H]} H. A. Hamm: Lefschetz theorems for singular varieties, Proc. Symp. Pure Math., Singularities, Volume 40, Part 1 (1983), 547-557.

\item{[L]} K. Lamotke: The topology of complex projective varieties
after S. Lefschetz, Topology 20(1981),15-51.

\item {[L\^e]} D.T. L\^e: Calcul du nombre de cycles \'evanouissants d'une hypersurface complexe, Ann. Inst. Fourier, Grenoble 23 (1973), 261-270.

\item{[LT]} D.T. L\^e, B. Teissier: Vari\'et\'es polaires locales et
classes de Chern des vari\'et\'es singuli\`eres, Ann. Math. 114 (1981),457-491.

\item{[N1]} A. N\'emethi: Th\'eorie de Lefschetz pour les vari\'et\'es
alg\'ebriques affines, C. R. Acad. Sci. Paris 303(1986), 567-570.

\item{[N2]} A. N\'emethi: Lefschetz theory for complex affine
varieties, Rev. Roum. Math. Pures et Appl. 33 (1988),233-260.

\item{[OT1]} P. Orlik, H. Terao: Arrangements of hyperplanes, Springer 1992.

\item{[OT2]} P. Orlik, H. Terao: Arrangements and Milnor fibers, Math. Ann. 301
(1995), 211-235.

\item {[OT3]} P. Orlik, H. Terao: The number of critical points of a product of powers of linear functions, Invent. Math. 120 (1995), 1-14.

\item{[PS]} S. Papadima, A. Suciu: Higher homotopy groups of complements of complex hyperplane arrangements, math.AT/0002251.

\item{[R]} R. Randell: Morse theory, Milnor fibers and hyperplane arrangements, math.AT 0011101.

\item{[ST]} D. Siersma, M. Tib\u ar: Singularities at infinity and their vanishing cycles II, Monodromy, Publ. RIMS, to appear.

\item{[T]} M. Tib\u ar: Asymptotic equisingularity and topology of complex hypersurfaces, Int. Math. Res. Not. 18 (1998), 979-990.

\bigskip

Laboratoire de Math\'ematiques Pures de Bordeaux

Universit\'e Bordeaux I

33405 Talence Cedex, FRANCE

\bigskip

email: dimca@math.u-bordeaux.fr

\bye